\documentclass[prd,eqsecnum,showpacs,nofootinbib]{revtex4}
\usepackage{amsmath}
\usepackage{bm}
\usepackage{amssymb}
\usepackage{graphicx}
\usepackage{epstopdf}

\begin{document}

\title{A Realization of ThurstonÕs Geometrization: Discrete Ricci Flow with Surgery\footnote{We dedicate this paper in honor of the 80th birthday of David Mumford for his leadership on the subject of computer vision.}}

\author{Paul M. Alsing$^{1}$,  Warner A. Miller$^{2}$\footnote{Corresponding Author: wam@fau.edu} \& Shing-Tung Yau$^{3}$}
\affiliation{
$^1$ Air Force Research Laboratory, Information Directorate, Rome, NY 13441\\
$^2$ Department of Physics, Florida Atlantic University, Boca Raton, FL 33431 \\
$^3$ Department of Mathematics, Harvard University, Cambridge, MA 02138
}

\begin{abstract}
Hamilton's Ricci flow (RF) equations were recently
  expressed in terms of a sparsely-coupled system of autonomous
  first-order nonlinear differential equations for the edge lengths
  of a $d$-dimensional piecewise linear (PL) simplicial geometry.
  More recently, this system of discrete Ricci flow (DRF) equations
  was further simplified by explicitly constructing the Forman-Ricci
  tensor associated to each edge, thereby diagonalizing the
  first-order differential operator and avoiding the need to invert
  large sparse matrices at each time step. We recently showed
  analytically and numerically that these equations converge for
  axisymmetric 3-geometries to the corresponding continuum RF
  equations.  We demonstrate here that these DRF equations yield an
  explicit numerical realization of Thurston's geometrization
  procedure for a discrete 3D axially-symmetric neckpinch geometry by
  using surgery to explicitly integrate through its Type-1 neckpinch
  singularity.  A cubic-spline-based adaptive mesh was required to
  complete the evolution.  Our numerically efficient simulations yield
  the expected Thurston decomposition of the sufficiently pinched
  axially symmetric geometry into its unique geometric structure --- a
  direct product of two lobes, each collapsing toward a 3-sphere
  geometry.  The structure of our curvature may be used to better
  inform one of the vertex and edge weighting factors that appear in
  the Forman's expression of Ricci curvature on graphs. 
  \end{abstract}

\maketitle

\section{Ricci Flow in 3D and its Applications}\label{sec:RFA}

Hamilton's Ricci flow (RF) yields new insights into a broad range of
problems from Perelman's proof of the Poincar\'e conjecture to
greedy-routing problems in cell phone networks
\cite{Hamilton:1982,Cao:2003,Chow:2004,Chow:2006,Chow:2007}. Here the
time evolution of the metric $\dot {\pmb g}$ is proportional to the
Ricci tensor ${\pmb Rc}$,
\begin{equation}
\label{eq:RF}
\dot {\pmb{g}} = -2\, Rc\left(\pmb{g}\right).
\end{equation}
The RF equation yields a forced diffusion equation for the curvature;
i.e., the scalar curvature ($R$) evolves as
\begin{equation}\label{eq:fdeqn}
\dot R = \triangle R + 2 R^2,
\end{equation}
here  $\triangle$ is the Laplacian with respect to the metric $\pmb{g}$.

The majority of the engineering applications of RF have been limited to the
numerical evolution of piecewise linear surfaces \cite{Gu:2012}.  This
is not surprising since a geometry with complex topology is most
naturally represented in a coordinate-free way by unstructured meshes,
e.g. finite volume \cite{Peiro:2005}, finite element
\cite{Humphries:1997}.  The applicability of discrete RF in two
dimensions arise from its diffusive curvature properties and from the
uniformization theorem for surfaces. Every simply connected
Riemann surface evolves under RF to one of three constant curvature
surfaces --- a sphere, a Euclidean plane or a hyperbolic plane.  RF on
surfaces is an accepted method for engineering a metric for a
surface given only its curvature \cite{Gu:2012}.  However, in three
dimensions, it is significantly more complicated.  In particular,
singularities can form during evolution under RF.  In three dimensions
the uniformization theorem yields the geometrization theorem of
Thurston showing that each closed 3-manifold has a similar, but richer,  
decomposition into a connected sum of one or more of eight
prime 3-manifolds \cite{Thurston:1997,Perelman:2003}.  The diffusive
curvature flow in three and higher dimensions together with this
classification provides a richer taxonomy than its 2-dimensional
counterpart.  We believe this more refined taxonomy may prove useful
in network classification.  Diffusive curvature flow can provide noise
reduction in higher dimensional manifolds, and in this direction we
are currently exploring a coupling of RF with persistent homology
\cite{Corne:2016}.  Finally, the soliton solutions of RF are Ricci
flat and are therefore vacuum solutions of Einstein's equations for
gravitation. This feature and its connection to the renormalization
group make RF with boundary an exciting topic for current research
into AdS/CFT models of quantum gravity \cite{ADS1,ADS2}.

\section{Discrete Ricci Flow in 3D}

\begin{figure}[ht] 
\centering
\includegraphics[width=5in]{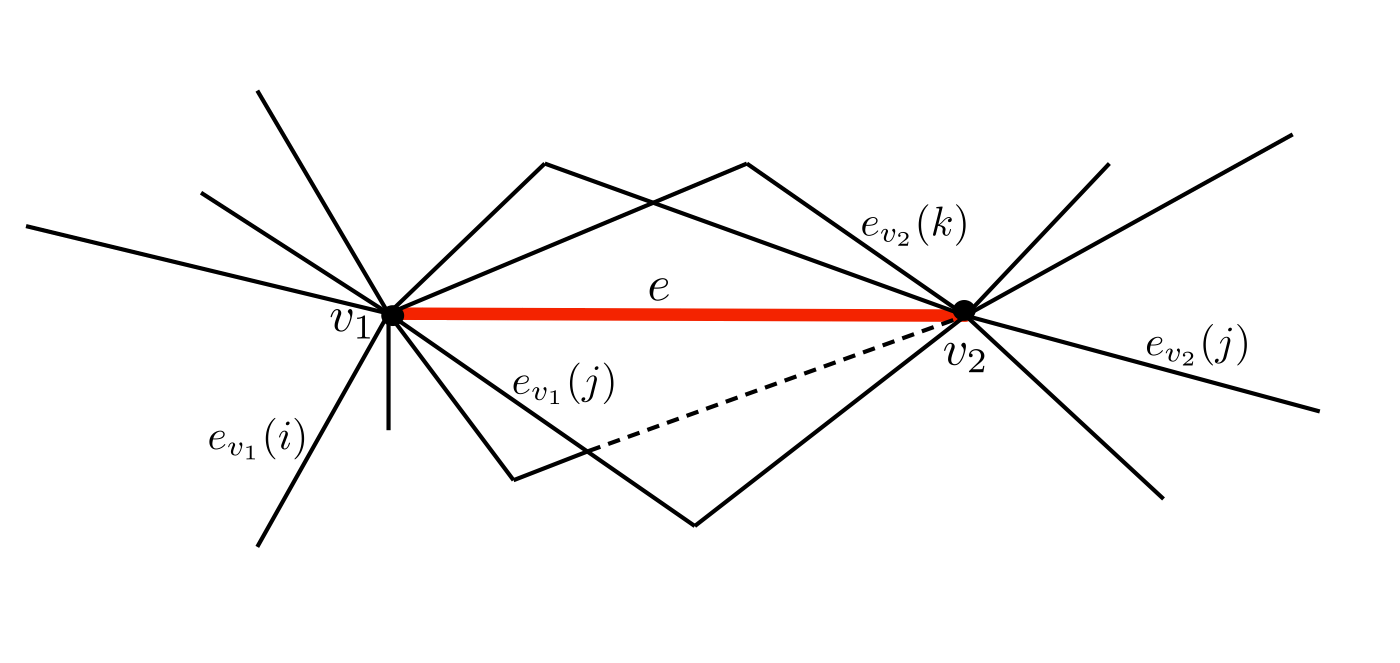} 
\caption{ The data structure of an edge ($e=\overline{v_1v_2}$) used
  in our definition of the Forman-Ricci tensor.  Included in this data
  structure are all the edges that share either vertex $v_1$ or $v_2$
  or both. This data structure is common for both the discrete Ricci
  flow tensor $Rc_{\!{}_{DRF}}$ used here as well as the Forman
  graph Ricci curvature $Rc_F$ for an edge $e$ in a graph described in Sec.~\ref{sec:higherD}.}
\label{fig:Kmv}
\end{figure}

A discrete RF (DRF) approach for three and higher dimensions, referred
to as Simplicial Ricci Flow (SRF), has been introduced recently and is
founded on Regge calculus \cite{Miller:2013,AMM:2011,McDonald:2012},
as well as complementary work in this direction by
\cite{Glickenstein:2011a,Glickenstein:2011,G:2005,Ge:2013,Forman:2003}.
The equations of SRF are similar to their continuum counterpart and
were shown for this model to convergent.  Recently, the $Rc_e$ tensor
was reconstructed on each edge $e$ of a lattice geometry
\cite{CM:2016} adhering closely to the approach by Forman
\cite{Forman:2003}. These new DRF equations form a diagonalized set of
first-order autonomous nonlinear differential equations in time. It is
numerically efficient and highly parallelizable. There is one equation
per edge in the lattice geometry,
\begin{equation}
\frac{1}{\ell_e} \frac{d\ell_e}{dt} = -Rc_e = -K_e + \frac{1}{2} R_e.
\label{DRF}
\end{equation}
In this DRF equation \cite{CM:2016}:
\begin{enumerate} 
\item $K_e$ is the sectional curvature of edge $\ell_e=\overline{v_1
  v_2} $ and is given in terms of the sum over all the edges,
  $\ell_{e_j}$ that share a common vertex ($v_1$ and/or $v_2$) with
  edge $\ell_e$,
\[
K_e = \sum_{e_{v_1},e_{v_2}\sim e} \frac{1}{2}\left(\frac{\cos^2(\theta_{e_{v1}}) \epsilon_{e_{v_1}}}{A_{e_{v_1}}}+
\frac{\cos^2(\theta_{e_{v2}}) \epsilon_{e_{v_2}}}{A_{e_{v_2}}}
\right).
\]
The data structure for this sectional curvature is illustrated in
(Fig.~\ref{fig:Kmv}).  It is expressed in terms of the Voronoi areas
$A_j$ dual to the edges $\ell_j$, the deficit areas $\epsilon_j$ of
these edges used in Regge calculus \cite{Regge:1961}, as well as the
angle $\theta_j$ between edge $\ell_e$ and $\ell_j$. Additionally,
\item $R_e$ is the scalar curvature associated
to edge $\ell_e$, and it is expressed in terms of the average of the
scalar curvatures at each of the endpoints of edge
$\ell_e=\overline{v_1 v_2} $,
 \[
 R_e = \frac{1}{2} \left( R_{v_1} + R_{v_2}\right).
 \]
 \end{enumerate}
 The vertex-based scalar curvatures were introduced earlier in Regge
 calculus, and is a certain weighted sum of the
 curvatures of the edges meeting a given vertex \cite{MM:2008},
 \[
 R_v = \frac{1}{V_v} \sum_{e\sim v}  \ell_e \epsilon_e.
 \]
Here $V_v$ is the dual volume associated with vertex $v$, and $\ell_e$
is the length of the edge emanating from vertex $v$.
 
It is the aim of this paper to explore the behavior of these new
diagonalized DRF equations in 3-dimensions for a geometry with axial
symmetry, and to examine the development of a Type-1 neck pinch
singularity through the singularity using manifold surgery
techniques. Thus providing the first piecewise linear numerical
realization of Thurston's geometrization using manifold surgery.

\section{The 3D Neckpinch Model}

\begin{figure}[ht] 
   \centering
   \includegraphics[width=3in]{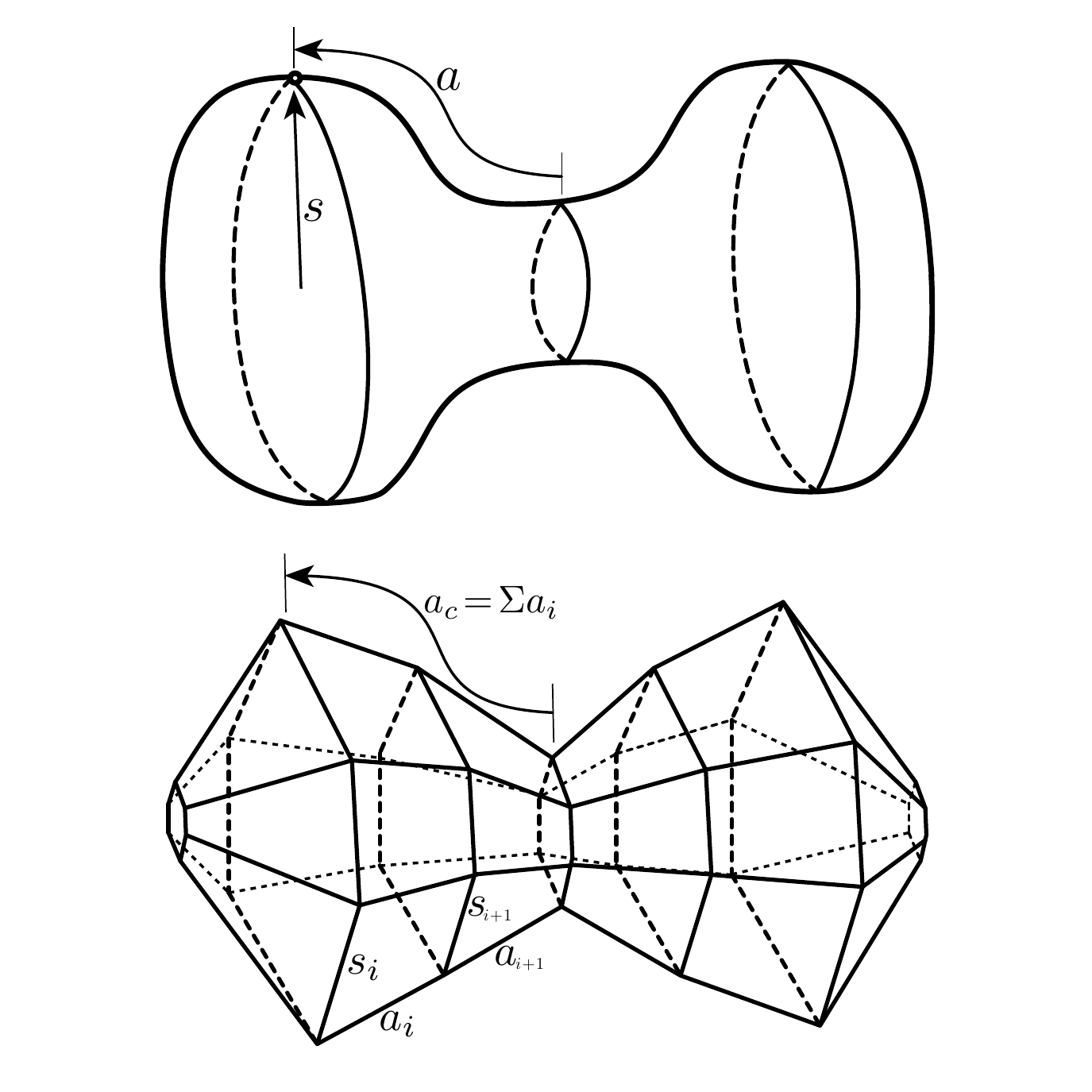} 
   \caption{A two dimensional representation of the 3D neckpinch
     geometry of Angenent and Knopf (continuum on top, and discrete on
     bottom).  In 3D the continuum cross-sections are 3-spheres and
     not circles, and in our discrete model the cross sections are
     icosahedrons and not hexagons.  The 3D cells are triangle-based
     frustum blocks as opposed to the trapezoids depicted in the
     bottom of the figure. Here the variable $a_c$ measures the proper
     distance from the equator, and $s$ is the length of the
     icosahedron edges.}
   \label{fig:DB}
\end{figure}

We use the analysis of Angenent and Knopf on the Type-1 singularity
analysis of the continuum RF equations as a foundation of this work
\cite{Knopf:2004}.  They carefully analyzed a class of axisymmetric
double-lobed shaped geometries with mirror symmetry about the plane of
the neck as illustrated in the top of Fig.~ \ref{fig:DB}.  The
symmetry of this geometry allows us to suppress one of the three
dimensions for visualization purposes.  In \cite{Knopf:2004} RF was 
applied to a warped product metric  on $I \times S^2$ having the form,
\begin{eqnarray}
\label{eq:wpm}
g & = &\underbrace{\varphi(z)^2 dz^2}_{da^2} + \rho(z)^2 g_{can}\\
 \label{eq:asmetric}
   & = & da^2 + \rho(a)^2 g_{can}.
\end{eqnarray}
Here, $I \in \mathbb{R}$ is an open interval,
\begin{equation}\label{eq:gcan}
g_{can} = d\theta^2 + \sin^2\theta d\phi^2,
\end{equation}
 is the metric of the unit 2-sphere,
\begin{equation}
a(z) = \int_{0}^{z} \, \varphi(z) dz,
\end{equation}
is the geodesic axial distance away from the waist, and $\rho(a)$ is
the radial profile of the mirror-symmetric geometry, i.e. $s=\rho(a)$
is the radius of the cross-sectional 2-sphere at axial distance $a$
from the waist.  Angenent and Knopf proved that the RF evolution
for such a geometry has the following properties:
\begin{enumerate}
\item If the scalar curvature is everywhere positive, $R\ge0$, then
  the radius of the waist ($s_{min}=\rho(0)$) is bounded, $(T-t) \le
  s_{min}^2 \le 2(T-t)$, where $T$ is the finite time at which a neck
  pinch occurs.
\item As a consequence, the neck pinch singularity occurs at or before
  $T = s_{min}^2$.
\item The heights of the two lobes are bounded from below and, under
  suitable conditions, the neck will pinch off before the lobes will
  collapse.
\item The neck approaches a cylindrical-type singularity.
\end{enumerate}
We demonstrated in our earlier work that the SRF equations, for a
sufficiently pinched radial profile, reproduced the neck pinch
singularity in finite time, and that the SRF evolutions agree with a
finite-difference solution of the continuum RF equations for the same
profile \cite{DB1:2014,DB2:2014}. However, in our previous analysis we
were unable to remove the singularity by manifold surgery and
so unable to integrate through the singularity and reproduce the direct product of
two collapsing 3-spheres. Furthermore the equations used previously,
though proven to converge to the continuum RF equations, form a
sparsely-coupled set of autonomous nonlinear first-order differential
equations that proved numerically difficult and time consuming to
solve.
\begin{figure}[ht] 
 \centering
   \includegraphics[width=5in]{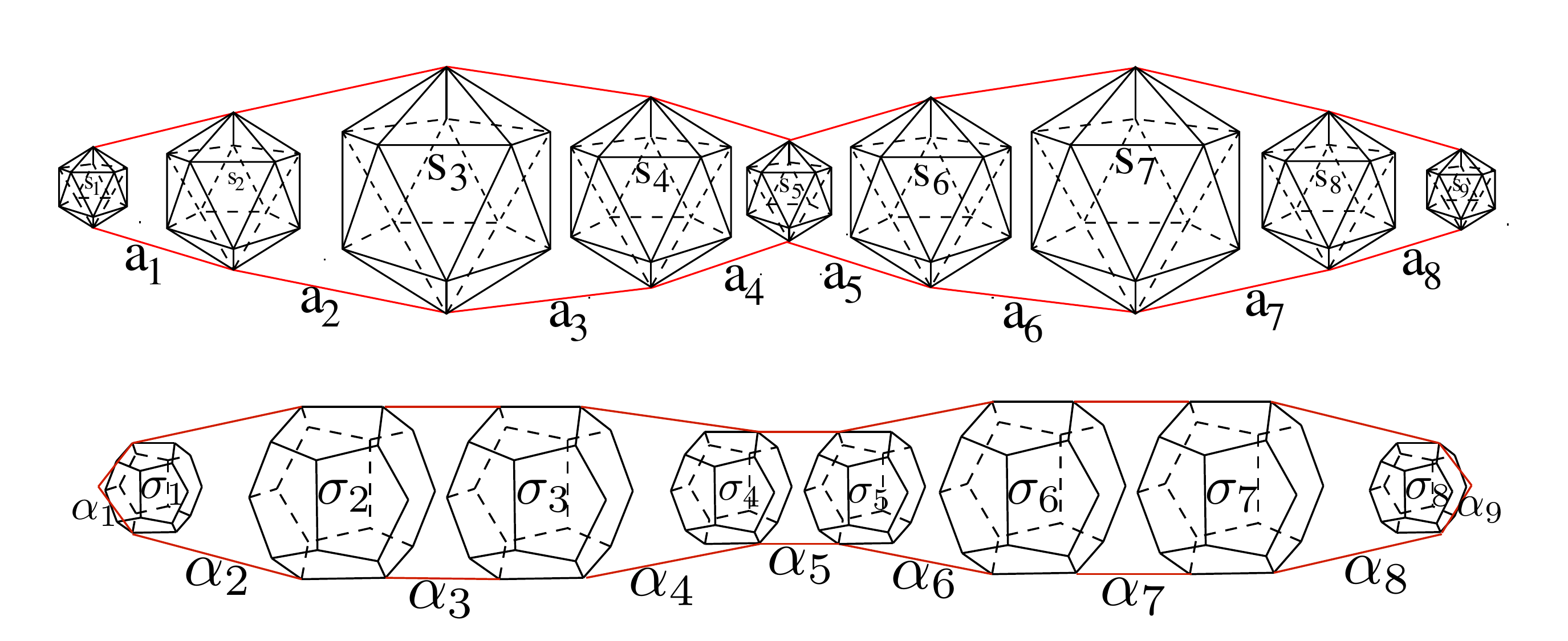} 
   \caption{An illustration of the icosahedron neckpinch geometry for
     nine cross-sectional icosahedra (top), and its dual dodecahedral
     lattice (bottom). The lattice is composed of triangle-based
     frustum blocks, and the dual lattice is composed of
     pentagonal-based frustum blocks.  The expressions for  the
     sectional, scalar, and Ricci curvature uses the dual lattice with its
     dodecahedral cross sections. }
   \label{fig:sigma}
\end{figure}

\begin{figure}[ht] 
\centering
\includegraphics[width=5.0in]{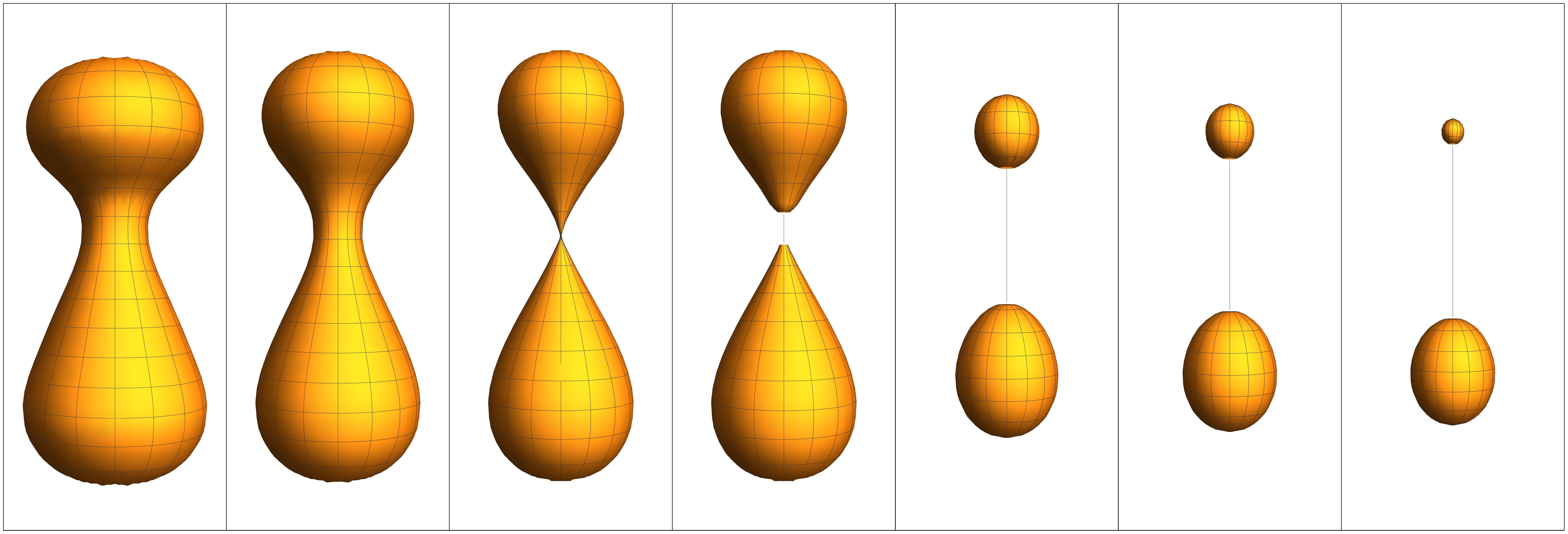} 
\caption{The RF of a lopsided neckpinch geometry through the Type-1
  singularity using surgery and yielding the geometry as a direct
  product of two 3-spheres.  We use axial symmetry of our model to
  suppress one dimension and the resulting two-lobed geometry can be
  visualized in Euclidean 3-space (our evolution was fortunately
  isometrically embeddable in $R^{3}$). The middle 3'rd and 4'th
  figure occur at the same time ($t=183.0$) in the evolution. They
  illustrate the explicit manifold surgery, where the spherical caps
  (two icosahedrons )are placed on the ends of the left and right
  lobes. This is the first numerical illustration of Thurston's
  geometrization procedure that we are aware of. This surface has 3438
  edges, 1580 triangle-based frustum blocks and 960 vertices, although
  symmetry reduces the number of edges to 80 icosahedral $\{s_i\}$
  edges and 79 axial $\{a_i\}$ edges.}
\label{fig:ellipse}
\end{figure}

The discrete model reported here is a piecewise linear (PL)
approximation to the double-lobed geometry (e.g. the $S^2$ cross
sections are modeled by icosahedra, and adjacent faces of the
icosahedra are connected to each other via frustum blocks) as
illustrated in Fig.~\ref{fig:sigma} and described more fully in
\cite{Miller:2013}.  Our simulation used 80 cross-sectional
icosahedra across the double-lobed profile.  We also relaxed the
condition of mirror symmetry about the throat and considered asymmetric
geometries.  This work represents the first non-trivial numerical
solution of the new DRF equations, and it is the first DRF integration
through a Type-1 singularity via manifold surgery of which we are aware.
The results are illustrated in Fig.~\ref{fig:ellipse}. Our algorithm is
numerically efficient, and the illustrative simulation presented here
involves the solution of a diagonal set of 159 autonomous nonlinear
first-order differential equations. We evolved the left and right lobes for 
1682 and 2133 time steps, respectively. We used a time step $\Delta t =0.25$. 
There is no longer the need for matrix inversion at each evolution step.

In the next section we describe the initial profile used and the numerical results
obtained.

\begin{figure}[ht] 
\centering
\includegraphics[width=3.5in]{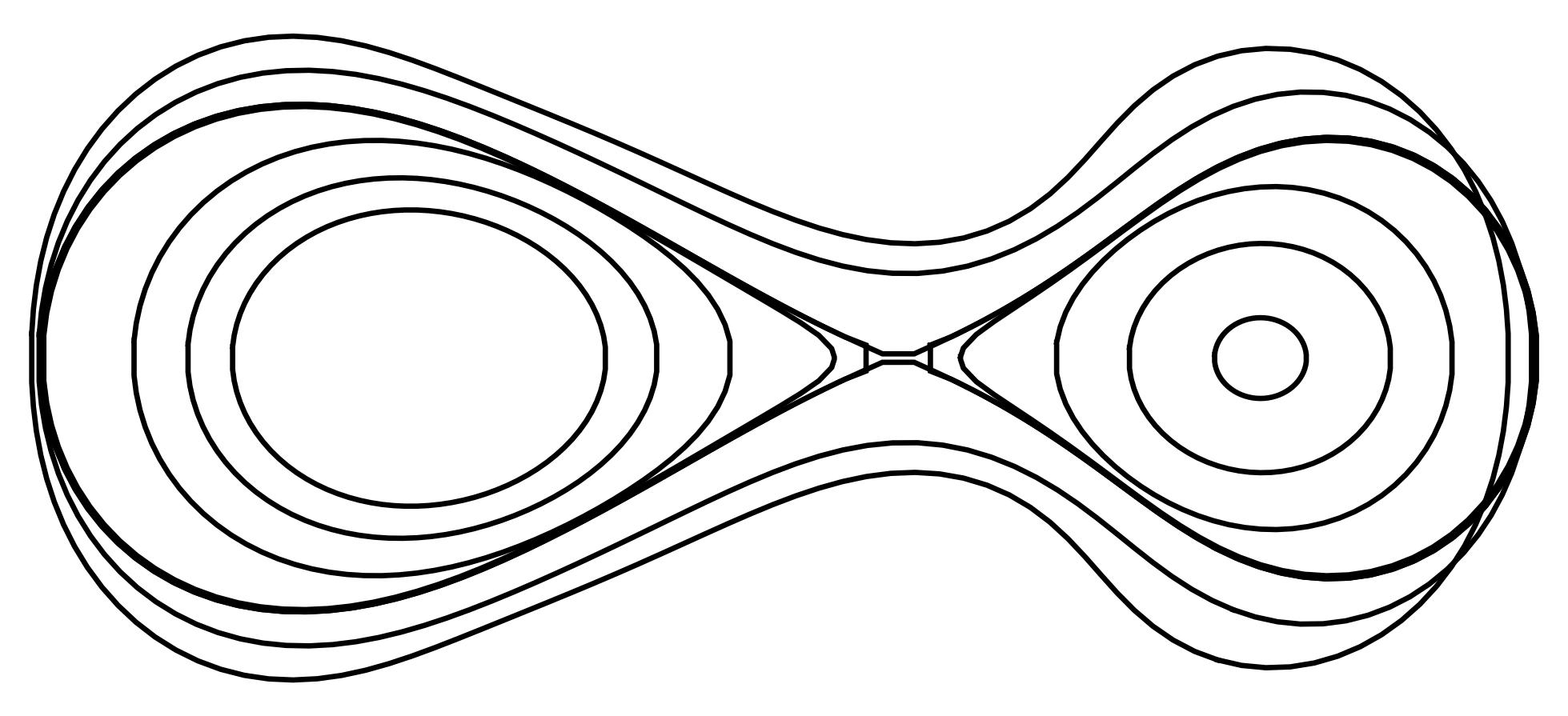} 
\caption{A 2-dimensional cross section of a lopsided neckpinch
  geometry evolving under RF through the Type-1 singularity.  Surgery
  yields two disconnected 3D ovoids and each becomes spherical under
  the RF evolution.  The resulting geometry is a direct product of two
  3-spheres. As the lobed geometry collapses a pinch occurs at t=
  183. At this point we remove the axial edges at the pinch and cap
  each end of the left and right lobe with a new icosahedra. These two
  surfaces (pre and post surgery) are the 3rd and 4th layers inside
  the initial surface. After surgery, we remesh both the left and
  right 3-dimensional ovoids using cubic spline interpolation.  This is, to our
  knowledge, the first numerical realization for PL manifolds of
  Thurston's geometrization procedure. This particular surface has
  3348 edges, 1580 triangle-based frustum blocks and 960 vertices,
  although symmetry reduces the number of edges to 80 icosahedral
  $\{s_i\}$ edges and 79 axial $\{a_i\}$ edges. }
\label{fig:ellipse2D}
\end{figure}



\section{DRF with Surgery: A Numerical Realization of Thurston's Geometrization for a Neckpinch Geometry.}
\label{sec:T}

\begin{figure}[ht] 
   \centering
   \includegraphics[width=2.75in]{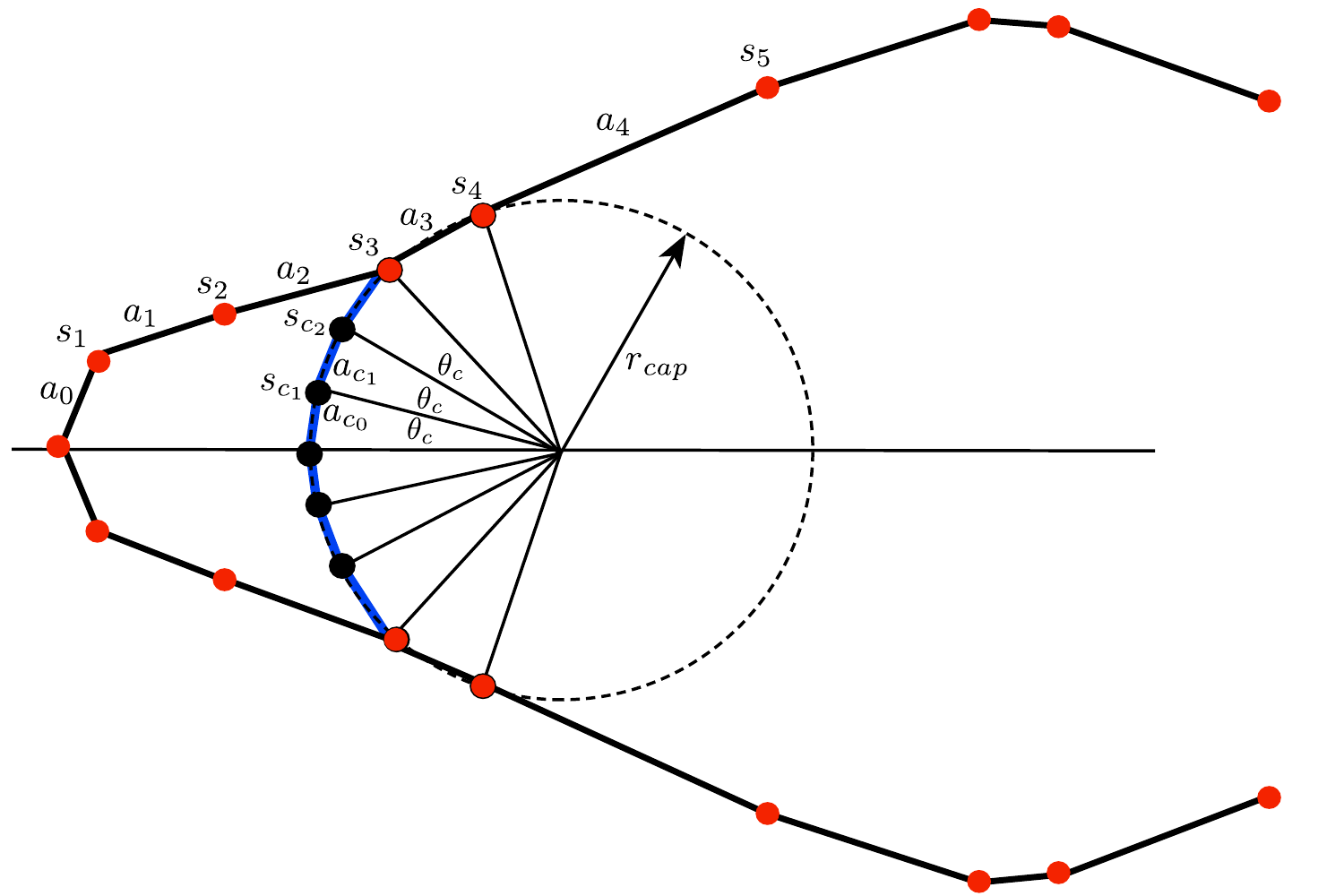} 
   \caption{ After the manifold surgery the lobe was closed using a
     spherical cap with proper matching conditions as illustrated in
     this figure. This involved reassigning the values to two of the
     $s$ variables and two of the $a$ values. This procedure offers no essential advantage
     over the simpler procedure consisting of just
     capping the surgery with an icosahedron and remeshing.  }
   \label{fig:cap}
\end{figure}

We evolved a sufficiently pinched axisymmetric 3-geometry which was
given the initial ($t=0$) radial profile,
\begin{equation}
s_i = 105.15 \left(1 - 0.2\  e^{{}^{\left(\frac{\xi_i + .4}{0.4}\right)^2}} - 0.05 \ e^{{}^{\left(\frac{\xi_i + 0.6}{0.3}\right)^2}} \cos(\xi_i) -  0.7 \cos(\xi_i)^4\right), \ \forall i\in\{1,...,n\},
\end{equation}
and axial segments,
\begin{equation}
a_i= 100\ \sin\left(\Delta \xi\right),\ \forall i\in\{1,2,...,n-1\},
\end{equation}
where $\xi_i = (n-2i+1)/2$, $\Delta \xi = \pi/(n+1)$, and there are $n=80$
icosahedral cross-sections.  
Fig.~\ref{fig:ellipse} shows the the initial profile of the lobed geometry in the rectangle to the left  
along with six other snapshots taken later during the evolution. This initial double-lobed
geometry is also illustrated in Fig.~\ref{fig:ellipse2D} and is the
outermost curve in the planar embedding. We evolve this surface by
numerically solving Eq.~\ref{DRF}.  This geometry evolved to a pinch
(third geometry from the left in Fig.~\ref{fig:ellipse}) at $t =
183.0$.  We evolved the equations using a fourth-order Runge-Kutta code
with $\Delta t=0.25$.  At every 50 steps in this evolution we remesh
the surface using a cubic spline interpolation.  This remeshing was
necessary to keep the circumcenter inside each frustum block (as
described in \cite{DB1:2014}).  Near the singularity $t=183$ we
removed the pinch by manifold surgery yielding the two lobes
exhibited in Fig.~\ref{fig:ellipse2D} using the following 4-step
procedure:
\begin{enumerate}
\item we remove the axial edge $a_{45}$ where the geometry pinched
  yielding a disconnected left and right lobed geometry each with
  $R^3$ topology (the right and left boundaries were removed;
  respectively);
\item we capped the left and right lobes by gluing an icosahedra to
  these open ends with edge length $s_{45}$ and $s_{46}$ thus forming
  two disconnected 3-dimensional ovoids;
\item we remeshed each of the 3-dimensional ovoids using a cubic
  spline; 
\item finally, we continued evolving using the DRF equations for both of the
  3-dimensional ovoids.
\end{enumerate}
A more sophisticated surgery procedure that we illustrate in
Fig.~\ref{fig:cap} was implemented.  Here we replace the last three
$s$ variables and two $a$ variables with their spherical cap values.
Because we found that this more time-consuming and sophisticated approach
yields the same results, we chose to use the more austere procedure
enumerated above. We evolved these two lobes separately using the
Eq.~\ref{DRF}. Under this flow the curvature uniformized and the lobes
each evolved toward a collapsing 3-sphere geometry as shown in the
figure.  We reproduced expected results with the new DRF equations as
shown in Fig.~\ref{fig:ellipse} and Fig.~\ref{fig:ellipse2D}. In other
words, the initial geometry evolved toward a direct product of two
constant curvature Thurston  geometries, and in particular, as a
direct product of two 3-spheres.

This numerical example demonstrates our ability to integrate
through a singularity and realize the Thurston decomposition.  It also
demonstrates that our current approach is numerically more efficient 
than our earlier formulations. 

\section{From Piecewise Linear Curvature to Graph Curvature}
\label{sec:higherD}

While, in this manuscript,  we have focused on the
discrete Ricci flow of a PL geometry and manifold surgery.  Our
formulation is based on Forman's curvature construction and can be
applied to more general structures, e.g. graphs.  It would be interesting to explore
the properties of graph curvature flow and determine its utility in
characterizing the graph structure, or in its ability to identify and
diffuse interesting curvature regions in the graph.  To this end,
there is considerable interest and pioneering work in applying the
Ricci flow techniques to characterize and identify change in dynamic
small-world spatial networks \cite{Weber:2016a,Weber:2016b}.  Positive
curvature networks stabilize, while negative hyperbolic curved
networks expand. The key to these approaches is a measure of the Ricci
curvature introduced by Forman \cite{Forman:2003}. We have identified
a striking, but intuitive, relationship between the Forman Ricci
curvature $Rc_F$ on graphs and our formulation of the discrete Ricci
tensor $Rc_{\!{}_{DRF}}$ \cite{Sree:2016,CM:2016},
\begin{align}
Rc_F &=  \frac{1}{2} \left( \frac{\omega(v_1)}{\omega(e)}+\frac{\omega(v_2)}{\omega(e)}\right)- \sum_{e_{v_1},e_{v_2}\sim e} \frac{1}{2} \left( \frac{\omega(v_1)}{\sqrt{\omega(e)\omega(e_{v_1})}}+\frac{\omega(v_2)}{\sqrt{\omega(e)\omega(e_{v_2})}} \right),\\
Rc_{\!{}_{DRF}} &= \frac{1}{2}\left(\frac{R_{v_1}+R_{v_2}}{2}\right) -\!\!\! \sum_{e_{v_1},e_{v_2}\sim e} \frac{1}{2}\left(\frac{\cos^2(\theta_{e_{v1}}) \epsilon_{e_{v_1}}}{A_{e_{v_1}}}+
\frac{\cos^2(\theta_{e_{v2}}) \epsilon_{e_{v_2}}}{A_{e_{v_2}}}
\right)
\end{align}
Here, $e$ is the edge under consideration between two nodes $v_1$ and
$v_2$, the edges sharing node $v_1$ are denoted by $e_{v_i}$ and are
each weighted by an appropriate weighting function $\omega(e_{v_i})\in
[0,1]$ (with $i=\{1,2\}$), and $\omega(v_i)\in [0,1]$ is the weighting
function for node $v_i$.  The data structure as shown in
Fig.\ref{fig:Kmv} is identical for both the discrete Ricci tensor and
the Forman curvature on graphs. The comparison of these two curvatures
for a given simplicial network, e.g. the 600-cell polytope, could
sharpen the definition of the vertex and edge weighting function for
the Forman curvature. This suggests the following correspondence:
\begin{align}
\omega(v_j) &\longleftrightarrow \cos^2\left(\theta_{e_j}\right) \epsilon_{e_j}\\
\sqrt{\omega(e_j)\omega(e)} &\longleftrightarrow A_{e_j}
\end{align}
We believe this may lead to discoveries characterizing
complex networks and work in this direction is already underway
\cite{Sree:2017}. 

It seems plausible that the set of DRF equations will have an
equally rich spectrum of application as does its 2-dimensional
counterpart known as combinatorial RF \cite{Chow:2003}.  We therefore
are motivated to explore the DRF in higher dimensions so that it can
be used in the analysis of topology and geometry, both numerically and
analytically, to bound Ricci curvature in discrete geometries and to
analyze and better handle higher--dimensional RF singularities
\cite{LinYau:2010,Knopf:2009}.  The topological taxonomy afforded by
RF is richer in 3D than in 2D.  In particular, the
uniformization theorem says that any 2--geometry will evolve under RF
to a constant curvature sphere, plane or hyperboloid, while in
3--dimensions the curvature and surface will diffuse into a connected
sum of eight distinct prime manifolds \cite{Thurston:1997}. We ask is
there a similar uniformization/geometrization theorem for 2D/3D
spatial networks?

\section*{Acknowledgements}
We wish to thank Rory Conboye and Matthew Corne for
  stimulating discussions and for their work.  We thank Rory Conboye for
  his help in reformulating the Forman-Ricci flow equations in their
  current numerically-efficient form.  PMA would like to acknowledge
  support of the Air Force Office of Scientific Research. We wish to
  thank the the Information Directorate of the Air Force Research
  Laboratory and the Griffiss Institute for providing us with an
  excellent environment for research. This work was supported in part
  through the VFRP and SFFP program, as well as AFRL grant
  \#FA8750-15-2-0047.  Any opinions, findings and conclusions or
  recommendations expressed in this material are those of the
  authors and do not necessarily reflect the views of AFRL.

\end{document}